\documentclass[psamsfonts, reqno, 11pt]{amsart}

\usepackage{amscd,amsmath,amssymb,amsfonts,verbatim}
\usepackage{times}
\usepackage[cmtip, all]{xy}
\usepackage{graphicx}
\usepackage[active]{srcltx}
\usepackage{color}
\usepackage{textcomp}

\definecolor{refkey}{gray}{.5}   
\definecolor{labelkey}{gray}{.5} 
\definecolor{Red}{rgb}{1,0,0}

\newcommand{\pf}{{\bf Proof : }}

\newcommand{\qedwhite}{\hfill \ensuremath{\Box}} 
\newcommand{\overbar}[1]{\mkern 1.5mu\overline{\mkern-1.5mu#1\mkern-1.5mu}\mkern 1.5mu}

\newtheorem{theo}{Theorem}[section]
\newtheorem{prop}[theo]{Proposition}
\newtheorem{lem}[theo]{Lemma}
\newtheorem{cor}[theo]{Corollary}

\theoremstyle{definition}

\newtheorem{notn}[theo]{Notation}
\newtheorem{rem}[theo]{Remark}
\newtheorem{exam}[theo]{Example}
\newtheorem{defi}[theo]{Definition}

\title{ Completion of skew completable unimodular rows }
\author{Sampat Sharma}
\newcommand{\Addresses}{{
  \bigskip
  \footnotesize

 \textsc{Sampat Sharma, Deaprtment of Mathematics, Indian Institute of Technology Bombay, \noindent
           Main Gate Rd, IIT Area, Powai, Mumbai, Maharashtra 400076, INDIA
   } \par\nopagebreak
  \textit{E-mail}: Sampat ~Sharma \texttt{<sampat.iiserm@gmail.com; sampat@math.iitb.ac.in>}

  \medskip

  }}

\begin{document}
\maketitle

\vskip0.15in

\subjclass 2010 Mathematics Subject Classification:{13C10, 13H99, 19B10, 19B14.}

~ {Local ring, regular local ring, Witt group}

\begin{abstract}
 In this paper, we prove that if $R$ is a local ring of dimension $d\geq 3$, $d$ odd and $\frac{1}{(d-1)!}\in R$ then any 
 skew completable unimodular row $v\in Um_{d}(R[X])$  is completable. It is also proved that skew completable unimodular rows
 of size $d\geq 3$ over a regular local ring of dimension $d$ are first row of a 2- stably elementary matrix.
\end{abstract}

\vskip0.50in

\begin{flushleft}
 Throughout this article we will assume $R$ to be a commutative noetherian ring with $1 \neq 0 .$
\end{flushleft}

\section{Introduction} In 1955, J.P. Serre
asked whether there were non-free projective modules over a polynomial
extension $k[X_1, \ldots, X_n]$, over a field $k$. D. Quillen (\cite{quill}) and A.A. Suslin (\cite{4}) 
settled this problem independently in early
1976; and is now known as the Quillen--Suslin theorem. Since every finitely generated projective module over $k[X_1, \ldots, X_n]$ is
stably free, to determine whether projective modules are free, it is enough to determine that unimodular rows 
over $k[X_1, \ldots, X_n]$ are completable. Therefore, problem of completion of unimodular rows is a central problem in 
classical $K$-Theory.

\par In \cite{swantowber}, R.G. Swan and J. Towber showed that if $(a^{2}, b, c)\in Um_{3}(R)$ then it 
can be completed to an invertible matrix over $R.$ This result of Swan and Towber was generalised by 
Suslin in \cite{sus2} who showed that if 
$(a_{0}^{r!}, a_{1}, \ldots, a_{r})\in Um_{r+1}(R)$ then it can be completed to an invertible matrix.
In  \cite{invent}, Ravi Rao studied the problem of completion of unimodular rows over $R[X]$,
where $R$ is a local ring. Ravi Rao showed that if 
$R$ is a local ring of dimension $d\geq2, \frac{1}{d!}\in R$, then any unimodular row over $R[X]$ of length
$d+1$ can be mapped to a factorial row by elementary transformations. In \cite{trans}, Ravi Rao proved that if $R$ is a local 
ring of dimension $3$ with $2R = R$, then unimodular rows of length $3$ are completable. In \cite{gr}, Ravi Rao generalised 
his result with Anuradha Garge and proved that if $R$ is a local ring of dimension $3$ with $2R = R$ then any unimodular 
row of length $3$ can be mapped to a factorial row via a two stably elementary matrix. 

\par In this article, we generalise the result of Garge--Rao for skew completable unimodular rows. We prove : 
\begin{theo}
 Let $R$ be a local ring of Krull dimension $d\geq 3$ with $d$ odd and $\frac{1}{(d-1)!} \in R.$ Let 
 $v = (v_{0}, v_{1},\dots,v_{d-1})\in Um_{d}(R[X])$ be skew-completable unimodular row over $R[X].$ Then there exists 
 $\rho \in SL_{d}(R[X])\cap E_{d+2}(R[X])$ and an invertible alternating matrix $W \in SL_{d+1}(R[X])$ such that 
 $$v\rho = e_{1}K(W).$$
\end{theo}

In the last section, we study the completion of skew completable unimodular rows over regular local rings. Since 
$SK_{1}(R[X])$ is trivial for a regular local ring $R$, we get the following result : 
\begin{theo}
 Let $R$ be a regular local ring of Krull dimension $d\geq 3$ with $d$ odd and $\frac{1}{(d-1)!} \in R.$ 
 Let $v = (v_{0}, v_{1},\dots,v_{d-1})\in Um_{d}(R[X])$ be skew-completable unimodular row over $R[X].$ Then there exists 
 $\rho \in SL_{d}(R[X])\cap E_{d+2}(R[X])$ such that $v = e_{1}\rho.$
\end{theo}

\section{Preliminary Remarks}
\par 
A row $v= (a_{0},a_{1},\ldots, a_{r})\in R^{r+1}$ is said to be unimodular if there is a $w = (b_{0},b_{1},\ldots, b_{r})
\in R^{r+1}$ with
$\langle v ,w\rangle = \Sigma_{i = 0}^{r} a_{i}b_{i} = 1$ and $Um_{r+1}(R)$ will denote
the set of unimodular rows (over $R$) of length 
$r+1$.
\par 
The group of elementary matrices is a subgroup of $GL_{r+1}(R)$, denoted by $E_{r+1}(R)$, and is generated by the matrices 
of the form $E_{ij}(\lambda) = I_{r+1} + \lambda e_{ij}$, where $\lambda \in R, ~i\neq j, ~1\leq i,j\leq r+1,~
e_{ij} \in M_{r+1}(R)$ whose $ij^{th}$ entry is $1$ and all other entries are zero. The elementary linear group 
$E_{r+1}(R)$ acts on the rows of length $r+1$ by right multiplication. Moreover, this action takes unimodular rows to unimodular
rows : $\frac{Um_{r+1}(R)}{E_{r+1}(R)}$ will denote set of orbits of this action; and we shall denote by $[v]$ the 
equivalence class of a row $v$ under this equivalence relation.

\subsection{The elementary symplectic Witt group $W_{E}(R)$} 
If $\alpha \in M_{r}(R), \beta \in M_{s}(R)$ are matrices 
then $\alpha \perp \beta$ denotes the matrix $\begin{bmatrix}
                 \alpha & 0\\
                 0 & \beta\\
                \end{bmatrix} \in M_{r+s}(R)$. $\psi_{1}$ will denote $\begin{bmatrix}
                 0& 1\\
                -1 & 0\\
                \end{bmatrix} \in E_{2}(\mathbb{Z})$, and $\psi_{r}$ is inductively defined by $\psi_{r} = \psi_{r-1}\perp 
                \psi_{1} \in E_{2r}(\mathbb{Z})$, for $r\geq 2$.
                \par 
    A skew-symmetric matrix whose diagonal elements are zero is called an alternating matrix. If $\phi \in M_{2r}(R)$ is 
    alternating then $\mbox{det}(\phi) = (\mbox{pf}(\phi))^{2}$ where $\mbox{pf}$ is a 
    polynomial (called the Pfaffian) in the matrix elements 
    with coefficients $\pm 1$. Note that we need to fix a sign in the choice of $\mbox{pf}$; so we insist 
    $\mbox{pf}(\psi_{r}) = 1$ 
    for all $r$. For any $\alpha \in M_{2r}(R)$ and any alternating matrix $\phi \in M_{2r}(R)$ we have 
    $\mbox{pf}(\alpha^{t}\phi
    \alpha) = \mbox{pf}(\phi)\mbox{det}(\alpha)$. For alternating matrices $\phi, \psi$ it is easy to check that 
    $\mbox{pf}(\phi \perp \psi) 
    = (\mbox{pf}(\phi))(\mbox{pf}(\psi))$.
    \par Two matrices $\alpha \in M_{2r}(R), \beta \in M_{2s}(R)$ are said to be equivalent (w.r.t. $E(R)$) if there exists  
    a matrix $\varepsilon \in SL_{2(r+s+l)}(R) \bigcap E(R)$,  such that $\alpha \perp \psi_{s+l} = \varepsilon^{t}
    (\beta \perp \psi_{r+l})\varepsilon,$ for some $l$. Denote this by $\alpha \overset{E}\sim \beta$. Thus 
    $\overset{E}\sim $ is an 
    equivalence relation; denote by $[\alpha]$ the orbit of $\alpha$ under this relation. 
    
    \par It is easy to see $(${\cite[p. 945]{7}}$)$ that $\perp$ induces the structure of an abelian group on the set of all 
    equivalence classes of alternating matrices with pfaffian $1$; this group is called elementary symplectic Witt group 
    and is denoted by $W_{E}(R)$.

\subsection{W. Van der Kallen's group structure on ${Um_{d+1}(R)}/{E_{d+1}(R)}$}
\begin{defi}{\bf{Essential dimension:}}
 Let $R$ be a ring whose maximal spectrum $\mbox{Max}(R)$ is a finite union of subsets $V_{i},$ where each $V_{i},$ 
 when endowed with the (topology induced from the) Zariski topology is a space of Krull dimension $d.$ We shall say 
 $R$ is essentially of dimension $d$ in such a case.
\end{defi}
\par 
For instance, 
a ring of Krull dimension $d$ is obviously essentially of dimension $\leq d;$ a local ring of dimension $d$ is essentially 
of dimension $0;$ whereas a polynomial extension $R[X]$ of a local ring $R$ of dimension $d\geq 1$ has dimension $d+1$ but is 
essentially of dimension $d$ as $\mbox{Max}(R[X]) = \mbox{Max}(R/(a)[X])\cup \mbox{Max}(R_{a}[X])$ for any non-zero divisor
$a\in R.$
\par
In $(${\cite[Theorem 3.6]{vdk1}}$),$ W. van der Kallen derives an abelian
 group structure on $\frac{Um_{d+1}(R)}{E_{d+1}(R)}$ when 
$R$ is essentially of dimension $d,$ for all $d\geq 2.$ Let $\ast$ denote the group multiplication henceforth. 
He also proved in $(${\cite[Theorem 3.16(iv)]{vdk1}}$),$ that the first 
row map is a group homomorphism
$$SL_{d+1}(R) \longrightarrow \frac{Um_{d+1}(R)}{E_{d+1}(R)}$$ 
when $R$ is essentially of dimension $d,$ for all $d\geq 2.$

\begin{lem}
\label{homo}
 Let $R$ be essentially of dimension $d\geq 2,$ and let $C_{d+1}(R)$ denote the set of all completable $(d+1)$-rows in 
 $Um_{d+1}(R).$ Then, 
 \begin{itemize}
  \item The map $\sigma \longrightarrow [e_{1}\sigma],$ where $e_{1} = (1,0,\ldots,0) \in Um_{d+1}(R),$ is a group homomorphism 
  $SL_{d+1}(R)\longrightarrow \frac{Um_{d+1}(R)}{E_{d+1}(R)}.$
  \item $\frac{C_{d+1}(R)}{E_{d+1}(R)}$ is a subgroup of $\frac{Um_{d+1}(R)}{E_{d+1}(R)}.$
 \end{itemize}
\end{lem}
  ${\pf}$ First follows from $(${\cite[Theorem 3.16(iv)]{vdk1}}$).$ Since $v\in C_{d+1}(R)$ can be completed to a matrix of 
  determinant one, $\frac{C_{d+1}(R)}{E_{d+1}(R)}$ is the image of $SL_{d+1}(R)$ under the above mentioned homomorphism; 
  whence is a subgroup 
  of $\frac{Um_{d+1}(R)}{E_{d+1}(R)}.$ 
  
  \begin{prop}
 Let $R$ be a local ring of dimension $d$, $d\geq 3$ and $\frac{1}{(d-1)!} \in R.$ Let $v = (v_{0},\ldots, v_{d}) 
 \in Um_{d+1}(R[X]).$ Then $v$ is completable if and only if $v^{(d-1)} = (v_{0}^{(d-1)}, v_{1},\ldots, v_{d})$ is completable.
\end{prop}
${\pf}$ In view of $(${\cite[Remark 1.4.3]{invent}}$),$ we may assume that $R$ is a 
reduced ring. By $(${\cite[Lemma 1.3.1, Example 1.5.3]{invent}}$),$ 
$$[v^{(d-1)}] = [v]\ast [v]\ast \cdots \ast[v],~~~~~~(d-1)~\mbox{times}$$
in $\frac{Um_{d+1}(R[X])}{E_{d+1}(R[X])}.$ By Lemma \ref{homo}, $v$ is completable implies $v^{(d-1)}$ is also completable.
\par Conversely, let $v^{(d-1)}$ be completable. By $(${\cite[Proposition 1.4.4]{invent}}$),$ 
$$v\overset{E}\sim (w_{0}, w_{1}, \ldots, w_{d-1}, c)$$
with $c\in R$ a non-zero-divisor. Since $\mbox{dim} (R/(c)) = d-1$ and $\frac{1}{(d-1)!} \in R,$ by 
$(${\cite[Corollary 2.3]{invent}}$),$
$$(\overset{-}w_{0}, \overset{-}w_{1}, \ldots, \overset{-}w_{d-1}) \in e_{1}SL_{d}(R/(c)[X]).$$ 
\par By $(${\cite[Proposition 1.2, Chapter 5]{mg}}$),$ $(w_{0}, w_{1}, \ldots, w_{d-1}, c^{d})$ is completable. Thus, 
\begin{itemize}
 \item $(v_{0}, v_{1}, \ldots, v_{d-1}, v_{d}^{d}) \overset{E}\sim (w_{0}, w_{1}, \ldots, w_{d-1}, c^{d})$ by 
 $(${\cite[Theorem 1]{orbit}}$),$
 \item $[v]^{n} = [(v_{0}, v_{1}, \ldots, v_{d-1}, v_{d}^{n})]$ for all $n$ by $(${\cite[Lemma 1.3.1]{invent}}$).$
\end{itemize}
\par 
Thus $[v]^{d} = [(w_{0}, w_{1}, \ldots, w_{d-1}, c^{d})] 
\in \frac{C_{d+1}(R[X])}{E_{d+1}(R[X])}$ and by hypothesis  $[v]^{d-1} = [v^{(d-1)}] \in \frac{C_{d+1}(R[X])}{E_{d+1}(R[X])}.$
Therefore by Lemma \ref{homo}, $v$ is completable.
$~~~~~~~~~~~~~~~~~~~~~~~~~~~~~~~~~~~~~~~~~~~~~~~~~~~~~~~~~~~~~~~~~~~~~~~~~~~~~~~~~~~~~~~~~~~~~~~~~~~~~~~~~~~~~~\qedwhite$

 \section{Krusemeyer's completion of the square of a skew completable row}
\setcounter {theo}{0}
 \begin{defi}
  A row $v\in Um_{2r-1}(R)$ is said to be skew completable if there is an invertible alternating 
  matrix $V \in GL_{2r}(R)$ with $e_{1}V = (0,v).$
 \end{defi}

First we note an example of skew completable unimodular row which is not completable.
\begin{exam}[Kaplansky]  Let $A = \frac{\mathbb{R}[x_{0}, x_{1}, x_{2}]}{(x_{0}^{2} + x_{1}^{2} + x_{2}^{2} - 1)}$ and $v = (\overline{x_{0}}, \overline{x_{1}}, \overline{x_{2}})\in Um_{3}(A).$ In view of $(${\cite[Section 5]{7}}$),$ every 
unimodular row of length 3 is skew completable. Thus $v = (\overline{x_{0}}, \overline{x_{1}}, \overline{x_{2}})$ is skew completable. Next we will show that $v$ is not completable.
\par Suppose to the contrary that $v = e_{1}\sigma$ for some $\sigma \in SL_{3}(A).$ Let $\sigma = (\sigma_{ij}).$ We can think $\sigma_{ij}$'s as a function on $S^{2}.$ Let us define tangent vector field
$$\phi : S^{2}\longrightarrow \mathbb{R}^{3}$$
$$~~~~~~~~~~~~~~~~~~~ w \longmapsto ((\sigma_{21}^{-1})^{t}(w), \sigma_{22}^{-1})^{t}(w), \sigma_{23}^{-1})^{t}(w)).$$ 
As $\sigma_{ij}$'s are polynomials, $\phi$ is a differential function. Since $(\sigma^{-1})^{t}\in SL_{3}(A)$, $\phi$ is a nonvanishing continuous tangent vector field on $S^{2}$ which is a contradiction to Hairy ball theorem. Thus $v$ is not completable.
\end{exam}
 \begin{theo}
 \label{kru}
  $($M. Krusemeyer$)$ $(${\cite[Theorem 2.1]{kruse}}$)$
  Let $R$ be a commutative ring and $v = (v_{1},\ldots, v_{n})$ be skew completable. Let $V$ be a skew completion of $v$, 
  then $(v_{1}^{2},v_{2},\ldots,v_{n})$ is completable.
 \end{theo}
\begin{notn} In the above theorem we will denote $K(V)\in SL_{n}(R)$ to be a completion of $(v_{1}^{2},v_{2},\ldots,v_{n})$ for a skew completable unimodular row $v = (v_{1},\ldots, v_{n})$ and its skew completion $V.$
\end{notn}
\begin{rem}
\label{kru1}
 M. Krusemeyer's proof in $(${\cite[Theorem 2.1]{kruse}}$),$ shows that $V \in (1\perp K(V)) E_{n+1}(R).$
\end{rem}
\begin{lem}
\label{kru2}
 Let $R$ be a commutative ring and $v = (v_{1},\ldots, v_{n})\in Um_{n}(R)$ be skew completable 
 to $V.$ Then $[e_{1}K(V)] = [e_{1}K(V)^{t}].$
\end{lem}
 ${\pf}$ By Remark \ref{kru1}, $V \in (1\perp K(V)) E_{n+1}(R).$
 Since $- I_{2k} \in E_{2k}(R),$ we have $V \in V^{t}E_{n+1}(R).$ Therefore $(1\perp K(V))^{t} \in (1\perp K(V)) E_{n+1}(R).$ 
 Since stably $K(V)$ and $K(V)^{t}$ are in same elementary class, therefore in view of $(${\cite[Lemma 10]{vas2}}$),$ 
 we have $[e_{1}K(V)] = [e_{1}K(V)^{t}].$
 $~~~~~~~~~~~~~~~~~~~~~~~~~~~~~~~~~~~~~~~~~~~~~~~~~~~~~~~~~~~~~~~~~~~~~~~~~~~~~~~~~~~~~~~~~~~~~~~~~~~~~~~~~
           ~~~~~\qedwhite$

 \begin{lem}
 Let $R$ be a local ring with $1/2\in R$ and let $V$ be an invertible alternating matrix of Pfaffian $1.$ 
 Let $e_{1}V = (0,v_{1}, \ldots, v_{2r-1}).
 $ Then $[V^{2^{n}}] = [W], ~\mbox{with} ~e_{1}W = (0,v_{1}^{2^{n}}, \ldots, v_{2r-1}).$
\end{lem}
${\pf}$ We will prove it by induction on $n.$ For $n = 1,$ by $(${\cite[Corollary 4.3]{gr}}$),$ 
$W_{E}(R[X])\hookrightarrow SK_{1}(R[X])$ is injective, we have 
$$V\perp V \overset{SK_{1}}\equiv V^{2} \overset{SK_{1}}\equiv V^{t}\psi_{r} V \overset{SK_{1}}\equiv
(1\perp K(V)^{t})\psi_{r} (1\perp
 K(V)).$$ 
 \par Therefore $[V^{2}] = [U]$ with $e_{1}U = (0,v_{1}^{2}, \ldots, v_{2r-1}).$ Now assume that result is true for all 
 $k\leq n-1$ and Let $[W_{1}] = [V^{2^{n-1}}]$ with $e_{1}W_{1} = (0,v_{1}^{2^{n-1}}, \ldots, v_{2r-1}).$ Since by lemma 
 $(${\cite[Corollary 4.3]{gr}}$),$  
 $W_{E}(R[X])\hookrightarrow SK_{1}(R[X])$ is injective, we have 
 $$W_{1}\perp W_{1} \overset{SK_{1}}\equiv W_{1}^{2} \overset{SK_{1}}\equiv W_{1}^{t}\psi_{r} W_{1} 
 \overset{SK_{1}}\equiv (1\perp K(W_{1})^{t})\psi_{r} (1\perp
 K(W_{1})).$$ 
  Therefore $[V^{2^{n}}] = [W]$ with $e_{1}W = (0,v_{1}^{2^{n}}, \ldots, v_{2r-1}).$
 $~~~~~~~~~~~~~~~~~~~~~~~~~~~~~~~~~~~~~~~~~~~~~~~~~~~~~~~~~~~~~~~~~~~~~~~~~~~~~~~~~~~~~~~~~~~~~~~~~~~~~~~~~
           ~~~~~\qedwhite$

\section{Completion of skew-completable unimodular rows of length ${d}$}
In this section, we prove that if $R$ is a local ring of dimension $d \geq 3$, $d$ odd and $\frac{1}{(d-1)!}\in R$ then any 
 skew completable unimodular row $v\in Um_{d}(R[X])$  is completable.
\begin{prop}
 \label{alt}
   Let $R$ be a local ring of dimension $d\geq 3$ with $d$ odd and $\frac{1}{(d-1)!} \in R.$ Let $V \in SL_{d+1}(R[X])$ be
   an 
   alternating matrix with Pfaffian $1.$ Then $[V] = [(1\perp K(W)^{t})\psi_{\frac{d+1}{2}}
   (1\perp K(W))]$ for some $[W]\in W_{E}(R[X]).$ 
   Consequently, there is a 1-stably elementary matrix $\gamma \in SL_{d+1}(R[X])$ such that $$V = \gamma^{t}
   (1\perp K(W)^{t})\psi_{\frac{d+1}{2}} (1\perp K(W))\gamma.$$
 \end{prop}
${\pf}$ By $(${\cite[Proposition 2.4.1]{trans}}$),$ $[V] = [W_{1}]^{2}$ for some $W_{1} \in W_{E}(R[X]).$ By 
$(${\cite[Theorem 2.6]{4}}$),$ $Um_{r}(R[X]) = e_{1}E_{r}(R[X])$ for $r\geq d+2,$ so on applying 
$(${\cite[Lemma 5.3 and Lemma 5.5]{7}}$),$ a few times, if necessary, we can find an alternating 
matrix $W\in SL_{d+1}(R[X])$ 
such that $[W_{1}] = [W].$ Therefore $[V] = [W]^{2}.$ Now, we have 
$$W\perp W \overset{SK_{1}}\equiv W^{2} \overset{SK_{1}}\equiv W^{t}\psi_{\frac{d+1}{2}} W \overset{SK_{1}}\equiv
(1\perp K(W)^{t})\psi_{\frac{d+1}{2}} (1\perp
 K(W)).$$  
 Since in view of $(${\cite[Corollary 4.3]{gr}}$),$ 
$W_{E}(R[X]) \hookrightarrow SK_{1}(R[X])$ is injective. Thus
$[W]^{2}= [(1\perp K(W)^{t})\psi_{\frac{d+1}{2}}
(1\perp K(W))].$ Therefore, $[V] = [(1\perp K(W)^{t})\psi_{\frac{d+1}{2}}
   (1\perp K(W))].$ The last statement follows 
by applying $(${\cite[Lemma 5.5 and Lemma 5.6]{7}}$).$
$~~~~~~~~~~~~~~~~~~~~~~~~~~~~~~~~~~~~~~~~~~~~~~~~~~~~~~~~~~~~~~~~~~~~~~~~~~~~~~~~~~~~~~~~~~~~~~~~~~~~~~~~~~~~~~\qedwhite$

\begin{theo}
\label{transgeneralise}
 Let $R$ be a local ring of Krull dimension $d\geq 3$ with $d$ odd and $\frac{1}{(d-1)!} \in R.$ Let 
 $v = (v_{0}, v_{1},\dots,v_{d-1})\in Um_{d}(R[X])$ be skew-completable unimodular row over $R[X].$ Then there exists 
 $\rho \in SL_{d}(R[X])\cap E_{d+2}(R[X])$ and an invertible alternating matrix $W \in SL_{d+1}(R[X])$ such that 
 $$v\rho = e_{1}K(W).$$
\end{theo}
${\pf}$ Let $V\in SL_{d+1}(R[X])$ be an invertible alternating matrix of Pfaffian $1$ which is a skew completion of 
$v = (v_{0}, v_{1},\dots,v_{d-1}).$ By Proposition \ref{alt}, there exists an alternating 
matrix $W\in SL_{d+1}(R[X])$ of 
Pfaffian $1$ such that 
$$[V] = [(1\perp K(W))^t\psi_{\frac{d+1}{2}} (1\perp K(W))].$$
Therefore there exists $\gamma \in SL_{d+1}(R[X])\cap E_{d+2}(R[X])$ such that $$\gamma^{t}V\gamma = 
(1\perp K(W))^t\psi_{\frac{d+1}{2}} (1\perp K(W)).$$
In view of $(${\cite[Corollary 5.17]{gr}}$),$ $e_{1}\gamma$ can be completed to an 
elementary matrix. Thus there exists $\varepsilon \in E_{d+1}(R[X]),$ and $\rho_{1} \in SL_{d}(R[X])\cap E_{d+2}(R[X])$ 
such that 
$$\varepsilon^{t}(1\perp\rho_{1})^{t}V(1\perp\rho_{1})\varepsilon = 
(1\perp K(W))^t\psi_{\frac{d+1}{2}} (1\perp K(W)).$$
By $($\cite[Corollary 4.5]{pr}$),$ there exists $\varepsilon_{1} \in E_{d}(R[X])$ such that 
$$(1\perp\varepsilon_{1})^{t}(1\perp\rho_{1})^{t}V(1\perp\rho_{1})(1\perp\varepsilon_{1}) = 
(1\perp K(W))^t\psi_{\frac{d+1}{2}} (1\perp K(W)).$$ Now we set $\rho = \rho_{1}\varepsilon_{1}.$
Thus $v\rho = e_{1}K(W)$. Hence $v$ is completable$.$
$~~~~~~~~~~~~~~~~~~~~~~~~~~~~~~~~~~~~~~~~~~~~~~~~~~~~~~~~~~~~~~~~~~~~~~~~~~~~~~~~~~~~~~~~~~~~~~~~~~~~~~~~~~~~\qedwhite$

\section{Completion of unimodular 3-vectors}
$(${\cite[Theorem 3.1]{trans}}$),$ Ravi A. Rao proved that for a local ring of dimension $3,$ every 
$v\in Um_{3}(R[X])$ is completable. We get stronger results than Theorem \ref{transgeneralise},
when we work with a local ring $R$ of dimension $3.$ We reprove Anuradha Garge and  
Ravi Rao's result in $(${\cite[Corollary 5.18] {gr}}$).$ 

\begin{prop}
 \label{altdim3} Let $R$ be a local ring of dimension $3$ with $\frac{1}{2k}\in R$ and let $V\in SL_{4}(R[X])$ be an
 alternating
  matrix of Pfaffian $1.$ Then $[V] = [V^{\ast}]$ in $W_{E}(R[X])$ with $e_{1}V^{\ast} = (0, a ^{2k}, b, c)$, 
  and $V^{\ast} \in SL_{4}(R[X]).$ Consequently, there is a stably elementary $\gamma \in SL_{4}(R[X])$ such that 
  $V = \gamma^{t}V^{\ast}\gamma.$ 
\end{prop}
${\pf}$  By $(${\cite[Proposition 2.4.1]{trans}}$),$ $[V] = [W_{1}]^{2k}$ for some $W_{1} \in W_{E}(R[X]).$ By 
$(${\cite[Theorem 2.6]{4}}$),$ $Um_{r}(R[X]) = e_{1}E_{r}(R[X])$ for $r\geq 5,$ so on applying 
$(${\cite[Lemma 5.3 and Lemma 5.5]{7}}$),$ a few times, if necessary, we can find an alternating 
matrix $V^{\ast}\in SL_{4}(R[X])$ 
such that $[W_{1}] = [V_{1}^{\ast}].$ Therefore $[V] = [V_{1}^{\ast}]^{2k}.$ Let $[V_{1}^{\ast}]^{2k} = [V^{\ast}],$ thus 
$[V] = [V^{\ast}].$
By $(${\cite[Lemma 4.8]{gr}}$),$ $e_{1}V^{\ast} = (0, a^{2k}, b, c).$ 
The last statement follows 
by applying $(${\cite[Lemma 5.3 and Lemma 5.5]{7}}$).$
$~~~~~~~~~~~~~~~~~~~~~~~~~~~~~~~~~~~~~~~~~~~~~~~~~~~~~~~~~~~~~~~~~~~~~~~~~~~~~~~~~~~~~~~~~~~~~~~~~~~~~~~~~~~~~~\qedwhite$

\begin{theo}
\label{transreprove}
 Let $R$ be a local ring of Krull dimension $3$ with $\frac{1}{2k} \in R.$ Let 
 $v = (v_{0}, v_{1}, v_{2}) \in Um_{3}(R[X])$. Then there exists 
 $\rho \in SL_{3}(R[X])\cap E_{5}(R[X])$  such that 
 $$v\rho = (a^{2k}, b,c )~\mbox{for~some}~(a,b,c)\in Um_{3}(R[X]).$$
\end{theo}
${\pf}$ Choose $w = (w_{0},w_{1},w_{2})$ such that $\Sigma_{i=0}^{2}v_{i}w_{i} = 1$, and consider the alternating matrix 
$V$ with Pfaffian $1$ given by 
$$V = \begin{bmatrix}
                 0 & v_{0} & v_{1} & v_{2}\\
                 -v_{0} & 0 & w_{2} & -w_{1}\\
                 -v_{1}& -w_{2} & 0 & w_{0}\\
                 -v_{2} & w_{1} & -w_{0} & 0\\
                \end{bmatrix} \in SL_{4}(R[X]).$$ By Proposition \ref{altdim3}, there exists an alternating 
matrix $V^{\ast}\in SL_{4}(R[X])$, with $e_{1}V^{\ast} = (0, a^{2k}, b, c),$ of Pfaffian $1$ such that 
$$[V] = [V^{\ast}].$$
Therefore there exists $\gamma \in SL_{4}(R[X])\cap E_{5}(R[X])$ such that $$\gamma^{t}V\gamma = V^{\ast}.$$
In view of $(${\cite[Corollary 5.17]{gr}}$),$ $e_{1}\gamma$ can be completed to an 
elementary matrix. Thus there exists $\varepsilon \in E_{4}(R[X]),$ and $\rho_{1} \in SL_{3}(R[X])\cap E_{5}(R[X])$ 
such that 
$$\varepsilon^{t}(1\perp\rho_{1})^{t}V(1\perp\rho_{1})\varepsilon = V^{\ast}.$$
By $($\cite[Corollary 4.5]{pr}$),$ there exists $\varepsilon_{1} \in E_{3}(R[X])$ such that 
$$(1\perp\varepsilon_{1})^{t}(1\perp\rho_{1})^{t}V(1\perp\rho_{1})(1\perp\varepsilon_{1}) = V^{\ast}.$$ 
Now we set $\rho = \rho_{1}\varepsilon_{1}.$
Thus $v\rho = (a^{2k}, b,c )~\mbox{for~some}~(a,b,c)\in Um_{3}(R[X])$.
$~~~~~~~~~~~~~~~~~~~~~~~~~~~~~~~~~~~~~~~~~~~~~~~~~~~~~~~~~~~~~~~~~~~~~~~~~~~~~~~~~~~~~~~~~~~~~~~~~~~~~~~~~~~~\qedwhite$

\section{Completion over regular local rings}
In this section, we prove that skew completable unimodular rows of size $d\geq 3$ over a regular local ring of dimension $d$ 
are first row of a 2- stably elementary matrix.

We note a result of Rao and Garge in $(${\cite[Corollary 4.3]{gr}}$).$ 
\begin{lem}
\label{gargerao}
 Let $R$ be a local ring with $2R = R.$ Then the natural map 
 $$W_{E}(R[X]) \hookrightarrow SK_{1}(R[X])$$ is injective.
\end{lem}

\begin{cor}
 \label{wittzero} Let $R$ be a regular local ring with $2R = R.$ Then the Witt group $W_{E}(R[X])$ is trivial.
\end{cor}
${\pf}$ Since $R$ is a regular local ring, $SK_{1}(R[X]) = 0.$ Thus the result follows in view of Lemma \ref{gargerao}.
$~~~~~~~~~~~~~~~~~~~~~~~~~~~~~~~~~~~~~~~~~~~~~~~~~~~~~~~~~~~~~~~~~~~~~~~~~~~~~~~~~~~~~~~~~~~~~~~~~~~~~~~~~~~~\qedwhite$

\begin{lem}
 \label{general} Let $R$ be a regular local ring of Krull dimension $d\geq 3$ with $d$ odd and $\frac{1}{(d-1)!} \in R.$ 
 Let $v = (v_{0}, v_{1},\dots,v_{d-1})\in Um_{d}(R[X])$ be skew-completable unimodular row over $R[X].$ Then there exists 
 $\rho \in SL_{d}(R[X])\cap E_{d+2}(R[X])$ such that $v = e_{1}\rho.$
\end{lem}
${\pf}$ Let $V$ be a skew completion of $v.$ In view of Corollary \ref{wittzero}, $W_{E}(R[X])=  0$, 
we have $[V] = [\psi_{\frac{d+1}{2}}].$ Thus upon applying $(${\cite[Lemma 5.5 and Lemma 5.6]{7}}$),$ there exists 
$\gamma \in SL_{d+1}(R[X]) \cap E_{d+2}(R[X])$ such that $\gamma^{t}V\gamma = \psi_{\frac{d+1}{2}}.$ 
\par By $(${\cite[Corollary 5.17]{gr}}$),$ $e_{1}\gamma$ can be completed to an 
elementary matrix. Thus there exists $\varepsilon \in E_{d+2}(R[X]),$ and $\rho_{1} \in SL_{d}(R[X])\cap E_{d+2}(R[X])$ 
such that 
$$\varepsilon^{t}(1\perp\rho_{1})^{t}V(1\perp\rho_{1})\varepsilon = \psi_{\frac{d+1}{2}}.$$
By $($\cite[Corollary 4.5]{pr}$),$ there exists $\varepsilon_{1} \in E_{d}(R[X])$ such that 
$$(1\perp\varepsilon_{1})^{t}(1\perp\rho_{1})^{t}V(1\perp\rho_{1})(1\perp\varepsilon_{1}) = \psi_{\frac{d+1}{2}}.$$ 
Now we set $\rho = (\rho_{1}\varepsilon_{1})^{-1}.$ Thus we have $v = e_{1}\rho.$
$~~~~~~~~~~~~~~~~~~~~~~~~~~~~~~~~~~~~~~~~~~~~~~~~~~~~~~~~~~~~~~~~~~~~~~~~~~~~~~~~~~~~~~~~~~~~~~~~~~~~~~~~~~~~\qedwhite$

\begin{cor}
 Let $R$ be a regular local ring of Krull dimension $3$ with $\frac{1}{2} \in R.$ Let 
 $v = (v_{0}, v_{1}, v_{2}) \in Um_{3}(R[X])$. Then there exists 
 $\rho \in SL_{3}(R[X])\cap E_{5}(R[X])$  such that $v= e_{1}\rho$.
\end{cor}
${\pf}$ Since every $v \in Um_{3}(R[X])$ is skew completable, thus the result follows in view of Lemma \ref{general}. 
$~~~~~~~~~~~~~~~~~~~~~~~~~~~~~~~~~~~~~~~~~~~~~~~~~~~~~~~~~~~~~~~~~~~~~~~~~~~~~~~~~~~~~~~~~~~~~~~~~~~~~~~~~~~~\qedwhite$

\Addresses


\begin{thebibliography}{9}





\bibitem{pr} 
P. Chattopadhyay, R.A. Rao;
\textit{Equality of elementary linear and symplectic orbits with respect to an alternating form}, 
J. Algebra 451 (2016), 46--64.

\bibitem{gr} 
A. S. Garge, R.A. Rao;
\textit{ A nice group structure on the orbit space of unimodular rows}, 
 K-Theory 38 (2008), no. 2, 113--133. 

 \bibitem{mg} 
S.K. Gupta, M.P. Murthy; 
\textit{Suslin's work on linear groups over polynomial rings and Serre problem}, no. 8, 
Indian Statistical Institute, New Delhi.



\bibitem{vdk1} 
W. van der Kallen;
\textit{A group structure on certain orbit sets of unimodular rows}, 
J. Algebra 82 (2) (1983), 363--397.



 \bibitem{kruse} 
M. Krusemeyer; 
\textit{Skewly completable rows and a theorem of Swan and Towber}, 
  Comm. Algebra 4 (1976), no. 7, 657--663.

 


\bibitem{quill} D. Quillen; \textit{Projective modules over polynomial rings},
Invent. Math. 36 (1976), 167--171.




\bibitem{invent}  R.A. Rao; \textit{ The Bass-Quillen conjecture in dimension three but
characteristic $\neq 2,3$ via a question of A. Suslin},
 Invent. Math. 93 (1988), no. 3, 609--618.

\bibitem{trans}  R.A. Rao; \textit{ On completing unimodular polynomial vectors of length three},
 Trans. Amer. Math. Soc. 325 (1991), no. 1, 231--239.






 


\bibitem{4} 
A.A. Suslin;
\textit{On the structure of special linear group over polynomial rings}, 
Math. USSR. Izv. 11 (1977), 221--238.

\bibitem{sus2} 
A.A. Suslin;
\textit{On stably free modules}, 
Math. USSR-Sb. 31 (1977), no. 4, 479--491.



\bibitem{swantowber} 
R.G. Swan, J. Towber;
\textit{A class of projective modules which are nearly free}, 
J. Algebra 36 (1975), 427--434.



\bibitem{7} 
L.N. Vaserstein, A.A. Suslin;
\textit{Serre's problem on projective modules}, 
Math. USSR Izv. 10 (1976), No. 5, 937--1001.



\bibitem{vas2} 
L.N. Vaserstein;
\textit{Computation of $K_{1}$ via Mennicke symbols}, 
 Comm. Algebra 15 (1987), no. 3, 611--656.
 
 
 \bibitem{orbit} 
L.N. Vaserstein;
\textit{operations on orbits of unimodular vectors}, 
 J. Algebra 100(2) (1986),  456--461.
 
 \bibitem{V} L.N. Vaserstein; \textit{On the stabilization of the general 
linear group over a ring},  Mat. Sbornik (N.S.) 79 (121) 405--424 (Russian); 
English translation in Math. USSR Sbornik. 8 (1969), 383--400.

\end{thebibliography}
           \end{document}